\documentclass{amsart}[12pt]
\usepackage[]{amsmath, amsthm, amsfonts, verbatim, amssymb}
\usepackage{tikz}
\usetikzlibrary{matrix,arrows}
\usepackage{url}
\usepackage[all,cmtip]{xy}

\usepackage{xcolor}

\usepackage{graphicx}
\usepackage{pinlabel}

\def\cE{{\mathcal E}}

\def\cH{\mathcal H}
\def\cO{{\mathcal O}}
\def\cI{{\mathcal I}}

\def\bZ{{\mathbb Z}}
\def\bQ{{\mathbb Q}}
\def\bC{{\mathbb C}}

\def\cS{{\mathcal S}}

\def\cC{{\mathcal C}}

\begin{document}
\newtheorem {theo}{Theorem}
\newtheorem {coro}{Corollary}
\newtheorem {lemm}{Lemma}
\newtheorem {rem}{Remark}
\newtheorem {defi}{Definition}
\newtheorem {ques}{Question}
\newtheorem {prop}{Proposition}
\def\spb{\smallpagebreak}
\def\mpb{\vskip 0.5truecm}
\def\bpb{\vskip 1truecm}
\def\wtM{\widetilde M}
\def\tM{\widetilde M}
\def\wtN{\widetilde N}
\def\tN{\widetilde N}
\def\tR{\widetilde R}
\def\tC{\widetilde C}
\def\tX{\widetilde X}
\def\tY{\widetilde Y}
\def\tE{\widetilde E}
\def\tH{\widetilde H}
\def\tL{\widetilde L}
\def\tQ{\widetilde Q}
\def\tS{\widetilde S}
\def\tc{\widetilde c}
\def\talpha{\widetilde\alpha}
\def\ti{\widetilde \iota}
\def\hM{\hat M}
\def\hq{\hat q}
\def\hR{\hat R}
\def\bs{\bigskip}
\def\ms{\medskip}
\def\ni{\noindent}
\def\td{\nabla}
\def\pd{\partial}
\def\hol{$\text{hol}\,$}
\def\Log{\mbox{Log}}
\def\bfQ{{\bf Q}}
\def\Todd{\mbox{Todd}}
\def\top{\mbox{top}}
\def\Pic{\mbox{Pic}}
\def\bP{{\mathbb P}}
\def\dxi{d x^i}
\def\dxj{d x^j}
\def\dyi{d y^i}
\def\dyj{d y^j}
\def\dzi{d z^I}
\def\dzj{d z^J}
\def\ozi{d{\overline z}^I}
\def\ozj{d{\overline z}^J}
\def\oz1{d{\overline z}^1}
\def\oz2{d{\overline z}^2}
\def\oz3{d{\overline z}^3}
\def\sI{\sqrt{-1}}
\def\hol{$\text{hol}\,$}
\def\ok{\overline k}
\def\ol{\overline l}
\def\oJ{\overline J}
\def\oT{\overline T}
\def\oS{\overline S}
\def\oV{\overline V}
\def\oW{\overline W}
\def\oY{\overline Y}
\def\oL{\overline L}
\def\oI{\overline I}
\def\oK{\overline K}
\def\oL{\overline L}
\def\oj{\overline j}
\def\oi{\overline i}
\def\ok{\overline k}
\def\oz{\overline z}
\def\om{\overline mu}
\def\on{\overline nu}
\def\oa{\overline \alpha}
\def\ob{\overline \beta}
\def\oGamma{\overline \Gamma}
\def\of{\overline f}
\def\oN{\overline N}
\def\og{\overline \gamma}
\def\ogamma{\overline \gamma}
\def\odelta{\overline \delta}
\def\otheta{\overline \theta}
\def\ophi{\overline \phi}
\def\opd{\overline \partial}
\def\oA{\overline A} 
\def\oB{\overline B}
\def\oC{\overline C}
\def\oD{\overline D}
\def\oIq1{\oI_1\cdots\oI_{q-1}}
\def\oIq2{\oI_1\cdots\oI_{q-2}}
\def\op{\overline \partial}
\def\ua{{\underline {a}}}
\def\us{{\underline {\sigma}}}
\def\Chow{{\mbox{Chow}}}
\def\vol{{\mbox{vol}}}
\def\dim{{\mbox{dim}}}
\def\rank{{\mbox{rank}}}
\def\diag{{\mbox{diag}}}
\def\tor{\mbox{tor}}
\def\supp{\mbox supp}
\def\bp{{\bf p}}
\def\bk{{\bf k}}
\def\a{{\alpha}}
\def\tchi{\widetilde{\chi}}
\def\ta{\widetilde{\alpha}}
\def\ovarphi{\overline \varphi}
\def\ocH{\overline{\cH}}
\def\tV{\widetilde{V}}
\def\tf{\widetilde{f}}
\def\th{\widetilde{h}}
\def\tT{\widetilde T}
\def\hG{\widehat{G}}
\def\hS{\widehat{S}}
\def\hD{\widehat{D}}
\def\Aut{\mbox{Aut}}
\def\hX{\widehat{X}}
\def\hC{\widehat{C}}
\def\hs{\widehat{s}}

\ni
\title[Bicanonical line bundle on compact complex $2$-ball quotients]
{Very ampleness of the bicanonical line bundle on compact complex $2$-ball quotients}
\author[Sai-Kee Yeung]
{Sai-Kee Yeung}

\begin{abstract}  
{\it The purpose of this note is to show that $2K$  of any smooth compact complex two
ball quotient is very ample, except possibly for four pairs of fake projective planes of minimal type, where $K$ is the canonical line bundle.
For the four pairs of fake projective planes, sections of $2K_M$ give an embedding of $M$ except possibly for at most two points on $M$.}
\end{abstract}

\address[]{Mathematics Department, Purdue University, West Lafayette, IN  47907USA} \email{yeung@math.purdue.edu
}

\thanks{\noindent{
The author was partially supported by a grant from the National Science Foundation}}

\ni{\it }
\maketitle

\begin{center}
{\bf 1. Introduction} 
\end{center}

\ms
\ni{\bf 1.1.}  General results on the very ampleness of the pluricanonical line bundle $pK_M$ for small $p$ on a complex manifold of general type
is a natural and interesting problem.
For a general smooth surface of general type, nice bounds have been obtained by Bombieri [Bom] and Reider [R].  The purpose of
this note is to investigate the smallest possible $p$ in the situation of smooth compact complex $2$-ball quotients, which 
are special but at the same time play prominent roles in some concrete algebraic and arithmetic problems.  In fact, they include as special cases
the Deligne-Mostow surfaces (cf. [P], [Le], [DM]), fake projective planes (cf. [M], [PY], [CS]) and Cartwright-Steger surfaces (cf. [CS]).
Geometrically, they form a boundary line in the geography of smooth surfaces of general type, sometimes known as the Bogomolov-Miyaoka-Yau line.
Arithmetically they contain one of the two main classes of Shimura varieties of dimension $2$, the other being the Hilbert modular surfaces.

By a compact complex $2$-ball quotient, we mean a smooth compact complex surface of the form $M=B_{\bC}^2/\Pi$, where
$B_{\bC}^2=\{(z_1,z_2)\vert |z_1|^2+|z_2|^2<1\}$ and $\Pi$ is a cocompact torsion-free lattice in $PU(2,1)$, the automorphism group of $B_{\bC}^2.$
It follows that the $M$ equipped with the Poincar\'e metric is a K\"ahler-Einstein metric of negative scalar curvature.  In particular, the canonical
line bundle $K_M$ is ample.  From Kodaira's Embedding Theorem, we know that the linear system $|pK_M|$ associated to the pluricanonical line
bundle $pK_M$ gives an embedding of $M$ if $p$ is sufficiently large.  In other words, $pK_M$ is very ample if $p$ is sufficiently large.  
The purpose of this note is to show the following sharper results in the case of complex $2$-ball quotients.  

\begin{theo}  Let $M$ be a smooth compact complex $2$-ball quotient.  Then $2K_M$ is very ample except possibly for fake projective planes
of minimal type.  For fake projective planes of minimal type, $2K_M$ is an embedding except possibly at at most two points.
\end{theo}

We recall that a fake projective plane is a complex surface with the same Betti numbers as but not biholomorphic to the complex projective plane.
A fake projective plane is a complex $2$-ball quotient, cf. [Y2].  It is said to be of {\it minimal type} in this article if the lattice involved is not contained properly as a 
sub-lattice in 
another lattice of $PU(2,1)$.  From the results of [PY] and [CS], there are only four pairs of fake projective planes of minimal type, given by
$(a=7,p=2,\{5\})$, $(a=7,p=2,\{5,7\})$, $(a=23,p=2,\emptyset)$ and $(a=23,p=2,\{23\})$ in the notation of Cartwright-Steger (cf. file \verb'registerofgps.txt' in the weblink of [CS]), see also
the table at the Appendix of this paper.

Theorem 1 is almost optimal since fake projective planes, which are complex $2$-ball quotients, have no sections in $K_M$.  
Furthermore, $K_M$ itself is asymptotically very ample on towers of coverings from [Y1].  We refer the readers to Lemma 1 in {\bf 2.1} for further elaborations.  Furthermore, the argument of \S10.5 of [PY] implies that $7H$ is very ample if $3H$ is numerically equivalent to $K_M$.

The result could be reduced to the following Theorem.  It corresponds to the special case of $M$ with Euler number $3$ and is
the marginal case in
the argument of Reider [R].

\begin{theo}
The bicanonical line bundle $2K_M$ is very ample for any smooth compact complex $2$-ball quotient $M$ with $c_2(M)=3$ apart from fake projective planes of minimal type.  
\end{theo}

The question for fake projective planes and the Cartwright-Steger surfaces is one of the original motivation for this note.  See also {\bf 2.1} for
other motivations.

Recall that sections of $qK_M$ can be regarded as automorphic forms of weight $q$ on a locally Hermitian symmetric space.  The result implies
that there are a lot of automorphic forms of weight $2$ on a compact complex $2$-ball quotient $M$.  In fact, they are enough to embed $M$ 
as a complex submanifold in some projective space apart from the four pairs of fake projective planes of minimal type.  

We also remark that there has been quite a lot of interesting work about birational properties of the bicanonical maps on a surface of general type.  
In particular, there is the result
of Borrelli in [Bor], which incorporates and summarizes earlier results of algebraic geometers including Catanese, Ciliberto, Debarre, Francia, Mendes Lopes, Morrison, Pardini, Xiao and others.
In particular, Mendes Lopes and Pardini proved in [MP] that $2K_M$ is a birational morphism for fake projective planes.  Very recently at the completion of the
first draft of the note, we were informed that
Di Brino and Di Cerbo had proved in [DD] the embedding of $2K_M$ for a fake projective plane with $|\Aut(M)|=21$.

For the embedding problem, the author does not know of a general method for smooth projective algebraic surfaces with $c_2(M)=3$. 
Apart from the argument of Reider in [R], we have to utilize the classification of fake projective planes and related surfaces in [PY], [CS] and [Y4], and 
exploit geometry of such surfaces, especially finite group actions.

The author is grateful to Lawrence Ein for raising the question about the Cartwright-Steger surfaces and for explaining the argument of Reider
to the author, to Rong Du and Ching-Jui Lai for helpful discussions, to Fabrizio Catanese and a referee for pointing out a gap in an earlier draft of this paper.  It is 
a pleasure for the author to express his gratitude to the referees for very helpful comments and suggestions on the article.

\bs
\begin{center}
{\bf 2. Preliminary discussions}
\end{center} 

\ms
\ni{\bf 2.1.}  To put the discussions in this note in perspective, let us collect some known facts.

\begin{lemm}
Let $M$ be a smooth compact complex $2$-ball quotient.\\
(a).  The Chern number $c_2(M)$ is a multiple of $3$ and $K_M^2$ is a multiple of $9$.\\
(b).  The divisor class $K_M$ is effective except for the case of $M$ being a fake projective plane.\\
(c). For any tower of coverings $\{M_i\}$ sitting above $M_1=M$, $K_{M_i}$ is 
very ample if $i$ is sufficiently large.

\end{lemm}
\ni{\bf Proof}  
(a). It is well-known that for a smooth compact complex $2$-ball quotient, a minimal surface of general type, the characteristic numbers
$c_1^2, c_2$ and $\chi(\cO)$ are all positive integers, cf. [BHPV].
 From Noether's Formula, we know that $\frac1{12}(c_1^2(M)+c_2(M))=\chi(\cO_M)\in \bZ^+$.
From the equality case of the Miyaoka-Yau inequality, we know that $c_1(M)^2=3c_2(M)$.
It follows that $c_2(M)\in 3\bZ^+$.  Hence $c_1^2(M)\in 9\bZ^+$.

(b). From definition, $1-h^1(M)+h^2(M)=\chi(\cO_M)\in \bZ^+.$  We see that  $h^2(M)=0$ if and only if
 $\chi(\cO_M)=1$ and $h^1(M)=0$, which implies that $c_2(M)=3$.  It follows that $M$ is a fake projective plane
as defined in [M] and classified in [PY].  We observe now that $h^0(M,K_M)=h^2(M)$ from Serre Duality.

(c).  By a tower of coverings of $M$, we mean a sequence of manifolds $M_i=M/\Pi_i$, where $\Pi_{i+1}$ is a normal subgroup of
$\Pi_{i}$ of finite index.  It is proved in [Y1] that $K_{M_i}$ is very ample if $i$ is sufficiently large.  We remark that it was shown earlier
by Hwang and To in [HT] that $2K_{M_i}$ is very ample if $i$ is sufficiently large.

\qed

The results prompt the question of what the smallest $p$ should be to make sure that $pK_M$ is very ample for all smooth compact
complex $2$-ball quotient $M$.  From Lemma 1b, we know that we need $p>1$ since $\Gamma(M,K_M)$ is trivial for $M$ being a
fake projective plane.  Hence Theorem 1 is the optimal result that we can expect.  For fake projective planes, it is known that $3K_M$ is
very ample, cf. [PY], {\bf 10.5}, as follows from a result of Bombieri [Bom].

\bs
\begin{center}
{\bf 3.   Tools from the results of Bogomolov and Reider.}
\end{center}

\ms
\ni{\bf 3.1.}  For the proof of Theorem 1, we recall first the following technical result, which is essentially a result of Reider [R], which is based on the result of Bogomolov (cf. [BHPV], p. 168).  We go through some details since the terminology and the argument are needed to prove our results in the next section.

\begin{prop}
Let $M$ be a smooth compact complex $2$-ball quotient for which with $c_2(M)=3$.  Then $2K_M$ is base point free.  Furthermore,
if sections of $\Gamma(M,2K_M)$ do not separate $p, q\in M,$ there exists an effective curve $B$ containing $p,q$, which can be two infinitesimally
close points, and satisfy
one of the following two properties,\\
(i). $B\cdot B=0$ and $K\cdot B=2$,\\
(ii). $K_M\equiv 3B$, where $`\equiv'$ stands for numerical equivalence.
\end{prop}

Before we go through the argument, we observe once again that $M$ is a smooth complex $2$-ball quotient with $c_2(M)=3$ and hence
$c_1^2(M)=9$.  From Noether Formula, $\chi(\cO_M)=1$.  Hence we conclude that the geometric genus $p_g(M)$ is the same
as the irregularity $q(M)$ of the surface.  Recall also that $K_M$ is ample since the Poincar\'e metric on $M$ has negative Ricci curvature.

For the embedding problem, in the following $Z$ is a subscheme of length $1$ or $2$ and would either be $Z_1=\{p\}$, $Z_2=\{p,q\}$, or $Z_3=\{2p\}$
corresponding to base point of $2K_M$, non-separated points or infinitesimally non-separable points.  Since the case of $Z=Z_1$ is already handled in
[R], we would focus on the cases of $Z=Z_2$ or $Z_3$.

As argued in [R], see also [BHPV] page 176, we may find a vector bundle $\cE$ coming from an extension
$$0\rightarrow \cO_M\rightarrow \cE\rightarrow \cI_Z\otimes K_M\rightarrow0 $$
with $c_1(\cE)=c_1(K_M)$ and $c_2(\cE)=\deg Z$.  Hence $c_2(\cE)= 2$ and 
\begin{equation}
c_1(\cE)^2-4c_2(\cE)\geqslant 1.
\end{equation}
Hence $\cE$ is Bogomolov unstable and we have from Bogomolov's Theorem (cf. [BHPV], p. 168) that there is a diagram

\bs
\begin{center}
\begin{tikzpicture}
implies/.style={double double equal sign distance, -implies},
\node (01) at (-4,0) {$0$};
\node (O1) at (-2,0) {$\cO_M$};
\node (E) at (0,0) {$\cE$};
\node (K) at (2,0) {$\cI_Z\otimes\cO(K_M)$};
\node (02) at (4,0) {$0$};
\node (v01) at (0,4) {$0$};
\node (L) at (0,2) {$\cO_M(L)$};
\node (B) at (0,-2) {$\cI_W\otimes\cO_M(B)$};
\node (v02) at (0,-4) {$0$};
\path[->] (01) edge node[above]{$$} (O1);
\path[->] (O1) edge node[above]{$$} (E);
\path[->] (E) edge node[above]{$$} (K);
\path[->] (K) edge node[above]{$$} (02);
\path[->] (v01) edge node[above]{$$} (L);
\path[->] (L) edge node[left]{$$} (E);
\path[->] (E) edge node[left]{$$} (B);
\path[->] (B) edge node[left]{$$} (v02);
\path[dashed,->] (O1) edge node[below]{$t$} (B);
\path[dashed,->] (L) edge node[above]{$t'$} (K);
\end{tikzpicture}
\end{center}
where $W$ is a $0$-dim scheme, and $L$ and $B$ are divisor line bundles satisfying
\begin{eqnarray}
L+B&=&c_1(\cE)=K_M\\
L\cdot B+\deg W &=&c_2(\cE)=\deg Z=2\\
(L-B)\cdot H&>&0 \ \forall \  \mbox{ample} \ H\\
(L-B)^2&>&4\deg W\geqslant 0
\end{eqnarray}

It is easy to see that $t$ is non-trivial and hence there is an effective divisor $B$ of $\Gamma(M,B)$ passing
through $Z$.   Let $d_1=K_M\cdot L$ and $d_2=K_M\cdot B$.  Then \begin{equation}
d_1+d_2=9, \ \ d_1>d_2>0.
\end{equation}
where the first identity follows from (2), and the second follows from (4).

Let $\delta=L\cdot B$.  From (2) and Hodge Index Theorem, we conclude that
\begin{equation}
\Delta:=d_1d_2-\delta(d_1+d_2)\leqslant 0.
\end{equation}
Hence $\delta>0$.  Identity (3) implies that $\delta\leqslant 2$.
Since
$$2p_a(B)-2=(K+B)\cdot B=(L+2B)\cdot B$$
in terms of arithmetic genus $p_a$ of $B$,  we conclude that $\delta$ is even
and hence $\delta=2$.  It follows from (3) that $W=0$ and hence $E$ is an extension of
two holomorphic line bundles.  Moreover, $L\cdot B=2$.  There are now two cases to consider.

\ms
\ni(i). $\Delta<0$.  In this case, recall from (7) that $d_1(d_2-2)-2d_2<0.$  Playing around with the integers satisfying (6),
we conclude that $d_2\leqslant 2$.  Suppose that $d_2=1$.  If follows $B\cdot B=K\cdot B-L\cdot B=-1.$  From 
the Adjunction Formula, $p_a(B)=1$, which contradicts hyperbolicity of $M$.  Hence we conclude that
$d_2=2$.  In such case $B\cdot B=0$ ad $p_a(B)=2$.

\ms
\ni(ii). $\Delta=0$.  The Hodge Index Theorem immediately implies that $L\equiv 2B$.  In such case, $K=L+B\equiv 3B$.



\qed

\begin{center}
{\bf 4. General argument  } 
\end{center} 

\ms
\ni{\bf 4.1.}   The main arguments of this article are in this and the next section.   The argument in this section is more general, and the next section
takes into account of classification and special geometric features.
We assume that $M$ is a smooth compact complex $2$-ball quotient with $c_2(M)=3$.

\begin{lemm}
Let $H$ be a line bundle on $M$ such that $3H=K \pmod\tor$ on $M$.  Let $\tau$ be any torsion line bundle.  Then $h^0(M,K+H+\tau)=3$.
\end{lemm}

\ni{\bf Proof}  From Riemann-Roch,
\begin{eqnarray*}
&&h^0(M,K+H+\tau)-h^1(M,K+H+\tau)+h^2(M,K+H+\tau)\\
&=&\frac12(K+H+\tau)\cdot(H+\tau)+\chi(M,\cO).
\end{eqnarray*}
Since $H$ is ample, $h^1(M,K+H+\tau)=0=h^2(M,K+H+\tau)$ from vanishing theorems.  As $\chi(M,\cO)=1$, 
it follows that
$h^0(M,K+H+\tau)=3.$

\qed

\ms
\ni{\bf 4.2.} 
\begin{lemm}
Assume that $3H=K+\sigma$ for some torsion line bundle $\sigma$.  Let $\tau$ be any torsion line bundle.  
Then $h^0(M,H+\tau)\leqslant 1$.
\end{lemm}

\ni{\bf Proof}  Assume that $H^0(M,H+\tau)$ contains two linearly independent sections $s_1$ and $s_2$.  Then
$s_1^4, s_1^3s_2, s_1^2s_2^2, s_1s_2^3, s_2^4$ give rise to five linearly independent sections of $K+H+\tau-\sigma$,
where $\tau-\sigma$ is torsion.  This follows from the fact that their generic vanishing orders along the zero divisor of $s_1$ are different.
This contradicts Lemma 2.

\qed

\ms
\ni{\bf 4.3.} 
We now prepare for the proof of Theorem 2.
 From the Chern number equality, we know that $K_M^2=9.$  It follows from the argument of [R] that $2K$ is base point free, as mentioned in Proposition 1.  Note from Riemann-Roch and Kodaira Vanishing Theorem
that $h^0(M,2K)=10$.   Hence $\Phi_{|2K_M|}:M\rightarrow N\subset P_{\bC}^9$
 is a morphism, where $N$ is the image of $\Phi_{|2K_M|}$.  
 
\ms 
We begin with the following lemma.

\begin{lemm} $N$ has dimension $2$. 
\end{lemm}

\ms
\ni{\bf Proof} Assume on the contrary that the image is a curve so that $\Phi_{|2K_M|}$
would give rise to a fibration.  From Stein Factorization, we may decompose $\Phi_{|2K_M|}=h\circ f$ where $f:M\rightarrow C$ has connected
fiber.    From Reider's
 result as in the proof of Proposition 1, each fiber  of the fibration would contain an effective divisor $B$ satisfying (i) or (ii) of Proposition 1.
 In case that $B$ is of type (ii), by considering a generic section,
 it follows that $B\cdot B= 0$, which contradicts the fact that $B\cdot B=\frac 19K_M\cdot K_M=1$.
 
 Hence it suffices to consider (i) with $B\cdot B=0$ and $K\cdot B=2$.  We have a fibration $f:M\rightarrow C$ to  a curve $C$ 
 with connected fibers of genus $2$ by  the Adjunction Formula.   Now as $M$ is a compact complex $2$-ball quotient, the fibration $f$ contains
 a singular fiber from  the result of [Li].  By considering base change and semi-stable reduction, we have a commutative diagram
 
 $$\begin{array}{ccc}
 N&\rightarrow&M\\
 \downarrow&&\downarrow\\
 C_1&\rightarrow&C
 \end{array}$$
 in which $f_1:N\rightarrow C_1$ is a semi-stable fibration with generic fiber of genus $2$.  An irreducible component of a singular fiber $R_1$ of $f_1$ is hence of genus $0$ or $1$.
 Hence from the commutative diagram and the Riemann-Hurwitz Formula, a singular fiber of $f$ contains a curve of genus $1$ or $0$, which leads to a contradiction since
 $M$ is hyperbolic.  The Lemma is proved.
 
 \qed

\ms
\ni{\bf 4.4.} It follows from the lemma 4 that $\Phi_{|2K_M|}$ is generically finite.  We claim the following result.
As mentioned in the Intoduction, for the special cases of fake projective planes, the result has already been proved earlier in 
 [MP].
 
\begin{lemm} $\Phi_{|2K_M|}$ is a birational morphism. 
\end{lemm}

 \ni{\bf Proof}    We already know that $\Phi_{|2K_M|}$ is base point free.
 Assume on the contrary that  $\Phi_{|2K_M|}$ is not birational.   We give two arguments.  The first invokes well established results in algebraic surfaces.
 Recall that a surface is called non-standard (cf. [Bor]) if $\Phi_{|2K_M|}$ is not birational and the surface contains no pencil of genus two curves.
 Hence $M$ is non-standard.

 
 From Theorem 0.7 of the result of Borrelli [Bor], we know that $\Phi_{|2K_M|}$ either has $p_g=q\leqslant 1$ or is a minimal model of a Du Val double plane, which means that 
 $M$ is birational to a twofold cover over $P_{\bC}^2$ or a Hirzebruch surface, cf. [BCP].  In the latter case, Theorem 0.4 of [Bor] implies that
 $M$ supports a rational pencil, which contradicts the fact that $M$ is hyperbolic.  Hence only the case that $p_g=q\leqslant 1$ occurs.
 
 In the case of $q=0$, $M$ is
 just a fake projective plane and there are only a finite choice of torsion line bundles  given by $H_1(M,\bZ)$ from Universal Coefficient Theorem.
 Since $y\in \Sigma$ is arbitrary, it implies that there are at least two different $y$ with the same $\tau_y$.  This implies that there are at least
 two linearly independent sections in  $\Gamma(M,L+\tau_y)$, contradicting Lemma 3.   
 
 In the case of $q=1$, the Picard variety of $M$ has
 complex dimension $1$.  On the other hand, the point $y\in \Sigma$, which has dimension $2$.  Hence there are at least two different $y_1, y_2\in \Sigma$
 for which $\tau_{y_1}=\tau_{y_2}$.  This implies that $\Gamma(M,L+\tau_{y_1})$ contains two sections $s_1, s_2$ passing through the two pairs
 of the inverse image of $y_1, y_2$ in $M$ respectively.  Again this contradicts Lemma 3 as above.

 The  argument above invoking [Bor] applies to much more general situation apart from complex $2$-ball quotients.

 An alternative and straight forward argument in our special  case of complex $2$-ball quotients is as follows.
 According to Proposition 1, either Case (i) or (ii) occurs.  Consider first Case (i).  
 The divisor $B$ in Proposition 1 has a moving part $B'$.  As $B^2=0$, we get
 $B'\cdot B'=0$ and $B'$ moves in base point free family $p:M\rightarrow C$.  Since $0<K_M\cdot B'\leqslant K_M\cdot B=2,$ we conclude that
 either $K_M\cdot B'=1$ or $2$.  The first case violates the Adjunction Formula, and for the second case, the Adjunction formula gives
 rise to genus $g(B')=2$.  As in the proof of Lemma 4, there exists a singular fiber in the family as $M$ is a complex $2$-ball
 quotient, which however leads to  a component of the singular fiber with normalization of genus $\leqslant 1$, contradicting the fact that $M$ is complex hyperbolic.
 We get a contradiction.  For Case (ii),  as $\Phi_{|2K_M|}$ is not birational, 
 we may still write $\Phi_{|2K_M|}:M\dasharrow\Sigma$ as a rational map of
  degree $d>1$.  Now Theorem 1(d) of [Y4]
  implies that for a smooth surface of general type $M$ with $c_2(M)=3$, $q=p_g$ and its value is either $0$ or $1$,  
  see also [Y5] for some details of the argument skipped in [Y4] and also some corrections.
 Once we know that $q=0$ or $1$, the arguments of the third and the fourth paragraphs of this section can be applied
 to reach a contradiction.

 \qed
 
\ms
\ni{\bf 4.5.} \begin{prop} Suppose $M$ is a smooth compact complex $2$-ball quotient with $c_1^2(M)=9$.  Then $\Phi_{|2K_M|}$
is an embedding except for a finite number of points lying on a fixed curve $B$.
\end{prop}

\ni{\bf Proof} 
We conclude from the result of Reider [R], as mentioned in Proposition 1, that $\Phi_{|2K_M|}$ is base point free.  From Lemma 5 we conclude that $\Phi_{|2K_M|}$ is a birational morphism.    
 We claim that $\Phi_{|2K_M|}$ is an immersion.  Let $D$ be a divisor on $M$ along which $\Phi_{|2K_M|}$ is not an immersion.  If 
 the dimension of $\Phi_{|2K_M|}(D)$ is $1$, the mapping $\Phi_{|2K_M|}$ is a branching map around $D$ of degree greater than $1$ and the argument of the last paragraph
 in the proof of Lemma 5 would
 rule it out, since we get a family of curves $B$ of type (ii) connecting preimages with respect to $\Phi_{|2K_M|}$ of points on the range.
  Suppose now that $D$ is an exceptional divisor with the image of  $\Phi_{|2K_M|}(D)$ being a point.  Grauert's Contraction Criterion
 (cf. [BHPV]) then implies that an irreducible component of $D_1$ satisfies $D_1\cdot D_1<0$.  However, according to (9) in Section 3, $D_1$ should
 be  a divisor  $B$ numerically the same as $\frac13K_M$ with positive self-intersection.  The contradiction implies that $\Phi_{|2K_M|}$ is an immersion
 and the claim is proved.  For embedding problem, the only case remained is that there may  be two different points mapped to the same
 image, which implies that the points have to lie on a curve $B$ as given in the case (ii) in Proposition 2.  From Lemma 3, such a curve $B$
is unique if exists.
 
 \qed

\begin{center}
{\bf 5. Geometry of specific ball quotients with $c_2=3$} 
\end{center} 

\ms
\ni{\bf 5.1.}  We conclude the proof of the main theorems in this section.  For this purpose, we utilize classification results in [PY], [CS]
and special geometric features of fake projective planes and the Cartwright-Steger surface.

Recall that a fake projective plane is an arithmetic complex $2$-ball quotient $\Gamma\backslash PU(2,1)/P(U(2)\times U(1))$.
We say that a fake projective plane or its associated lattice  $\Gamma$ is of {\it minimal type} if the ball quotient does not cover any other 
ball quotient apart from itself, or equivalently, $\Gamma$ is not properly contained in any other lattice $\Lambda$
of $PU(2,1)$.   A fake projective plane is said to be of {\it non-minimal type} if it is not of minimal type.  
From the classification of [PY] and [CS], there are only four lattices of {\it minimal type}.

\begin{theo}  
Let $M$ be a fake projective plane of non-minimal type or a Cartwright-Steger surface.  Then $2K_M$ is very ample.
\end{theo}
\ni {\bf Proof}  Lacking of a uniform proof, we have to consider separate arguments for different cases.  The idea is to make use of the specific geometry of the surfaces involved.  We are going to use extensively
the table of fake projective planes contained in [CS], the geometry of fake projective planes in [PY] and geometry of the Cartwright-Steger surface 
in [CS] and [CKY].  We separate the surfaces into different cases and prove the results one at a time in the following subsections.  A table is provided in the 
appendix which lists each fake projective plane according to the cases (b) to (d) below.

\ms
\ni{\bf 5.2.} {\it  Case (a). Cartwright-Steger surface:} 

\ms
Let $M$ be a Cartwright-Steger surface.  We refer the readers to [CS] for general facts about the lattice associated to the surface.  See also [CKY] for more
algebraic geometric properties.
Assume that $2K_M$ is
not very ample.  Proposition 2 implies that $K_M$ is linear equivalent to $3H$.  
From the explicit description of the lattice $\Pi$ in [CS], 
the abelianisation $H_1(M,\bZ)$ is found to be $\bZ\oplus\bZ$. Hence from the Universal Coefficient Theorem, $M$ has no torsion element apart from $\Pic^0(M)$. Hence 
we may write $K_M=3H+\tau$, where $\tau$ is an element in $\Pic^0(M)$,
which is a one torus. By taking division on the one torus, we know we can write $\tau=3\sigma$ for some torsion line bundle $\tau$. It follows that $K_M=3(H+\sigma)$. The argument of 10.4 of [PY] implies that the fundamental group $\Pi=\pi_1(M)$ can be lifted as a lattice from $PU(2,1)$ to $SU(2,1)$. This however contradicts a fact from [CS] that such lifting is not possible from explicit description of $\Pi$. Hence Proposition 2 implies that $2K_M$ is very ample.

\ms
\ni{\bf Remark} In the following, we consider fake projective planes.  Results in \S4 would imply that $2K_M$ is very ample for all fake projective planes if Conjecture 2 of [LY] is proved,
which is equivalent to $h^2(M,2B)=0$ in our notation.  Hence if the conjecture is true, the discussions in \S4, {\bf 5.2} and [Y4] would imply that $2K_M$ is very ample for
all smooth complex $2$-ball quotients.  At this point, the conjecture was proved only for cases with $|\Aut(M)|=9$ and $21$.  We refer the reader to [LY] for a proof and for
related references.

\ms
\ni{\bf 5.3.} {\it Case (b). Fake projective planes with a non-trivial automorphism group of order $3$:}

\ms

Let $M$ be a fake projective plane with $\Aut(M)\neq\{1\}$.  We know from the results of Keum [K] and Cartwright-Steger [CS] that $\Aut(M)$ 
is either trivial, or has order given by $21, 9$ or $3$.  Main Theorem' of [LY] and Proposition 2 immediately conclude the proof for 
fake projective planes with automorphism group of order $9$ and $21$, for the latter case, we use the fact that every torsion line bundle on $M$ is 
a $2$-torsion line bundles on $M$ according to the
file  \verb'registerofgps.txt' on the weblink of [CS].  Nonetheless the argument below gives an alternate proof for all these cases as well,
since $\Aut(M)$ contain a subgroup of order $3$.

Let $p, q$ be two points
which are not separated by $\Phi_{|2K_M|}$ and are contained in $B$.   Since the irregularity $q(M)=0$ from definition of a fake projective plane, Proposition 1(ii) can be written as $K_M=3B+\tau$ for some torsion line bundle $\tau$
where $K_M=3B+\tau$ and $\tau\neq0$.  $B\cdot B=1$ from definition.  From Lemma 6 of [LY], we know that $B$ is smooth of genus $3$.


Suppose now that $\Aut(M)=\bZ_3$ or more generally that $\Aut(M)$ contains a group of order $3$ denoted by $\bZ_3$.  Denote $X=M/\bZ_3$.  Let $\sigma$ be the generator of $\bZ_3$ and $p:M\rightarrow X$ be the
regular covering map.
We know that $B\in \Gamma(H+\eta)$ for some torsion line bundle $\eta$.   By going through the table of fake projective planes given by
the
file  \verb'registerofgps.txt' on the weblink of [CS], we see case by case that $H+\eta$ is invariant under $\bZ_3$ as a line bundle.  
Note that $\eta$ corresponds to an element in $H_1(M,\bZ)$.  If $\eta$ belongs to the subgroup $p^*H_1(X,\bZ)$, clearly $H+\eta$ is invariant under $\bZ_3$.
On the other hand, as $(\sigma^3)^*B=B$, we see that $\eta\in [H_1(M,\bZ)]^{\bZ_3}$, the part of the cohomology group invariant under $\bZ_3$.  By checking over each case in the list of fake projective planes,
we verify that $[H_1(M,\bZ)]^{\bZ_3}=p^*H_1(X,\bZ)$.  We conclude that $\eta$, and hence $H+\eta$, is invariant under $\bZ_3$.

It follows from Lemma 3 that $B$ is unique and is invariant as a set under $\Aut(M)$.   Let $C=B/\bZ_3$
and $g(X)$ be the geometric genus of a curve $X$.  It follows that $C$ is smooth, and from Riemann-Hurwitz's Formula,
\begin{equation}
2(g(B)-1)=3\cdot 2(g(C)-1)+\sum_{i=1}^l b_i
\end{equation}
where the sum is over fixed points of $\Aut(M)$ and $b_i=2$ is the ramification order, and $l\leqslant 3$ is the number of ramification points on
$B$.  Here $g(B)=3$ from Lemma 6 of [LY].  It follows that $g(C)=1$ and $k=2$.  Hence $C$ is actually a smooth
elliptic curve.  Moreover, $K_X$ is a $\bQ$-line bundle and $K_X\cdot C=\frac13 K_M\cdot B=1.$  

In cases that the hyperplane line bundle $H$ on $B_{\bC}^2$ descends to $X$, $C$ is numerically equivalent to the hyperplane $H$.
We observe from the discussions in [PY] and [CS] that apart from fake projective planes arising from the classes of $\cC_{18}$ in
the notation of [PY],  $H$ descends to $X$ as a line bundle, corresponding to the 
fact that $\Gamma$ can be lifted to $SU(2,1)$.  In the case of the classes of $\cC_{18}$, we still know from Poincar\'e Duality that
$C\equiv H'$ on $X$,  where $H'$ is a generator of the Neron-Severi group.  Here we used the fact that $h_2(M,\bZ)=1$.

Let $\pi:\hX\rightarrow X$ be the minimal desingularisation of $X$.  

\begin{lemm}
$X$ has three rational singularities of type $\frac13(1,2)$. Let $\tau:\hX\rightarrow X$ be the minimal resolution of $X$.  Then $K_{\hX}=\tau^*K_X$, $K_{\hX}^2=3$, and $h^0(X,2K_X)=4$.
\end{lemm}

\ni{\bf Proof}  This is essentially contained in Case 1 of \S4 in [Y3], which in turn depends on the information about singularity type of $\bZ_3$ action on a fake projective plane
given by Cartwright-Steger [CS] and Keum [K].  Hence it is known that $X$ has three rational singularities $Q_1, Q_2, Q_3$ of type $\frac13(1,2)$.  In the resolution, $\tau:\hX\rightarrow X$, each of 
$Q_i$ gives rise to a chain of two $(-2)$-curves $E_{ij}, j=1,2$.  One checks that $K_{\hX}=\tau^*K_X$ since the singularities are
rational, and that $K_{\hX}$ is nef and big from the fact that $X$ has
Picard number $1$.  It follows from Riemann-Roch and Kawamata-Viehweg Vanishing Theorem that 
$h^0(X,2K_X)=4$.

\qed

\ms
  Since $l=2$ in (8), we may assume from the above discussions
that $C$ passes through $P_1$ and $P_2$ but not $P_3$.  Let $\hC$ be the proper transform of
$C$ on $\hX$.  We may write 
\begin{equation}
\hC=\tau^*C-\sum_{i=1}^2\sum_{j=1}^2a_{ij}E_{ij}
\end{equation}
for some rational numbers $a_{ij}$.

Let the $\bZ_3$ action on a neighborhood of a singularity of type $\frac13(1,2)$ be given by 
$\sigma(x,y)=(\omega x, \omega^2y)$, where $\omega$ is a cubic root of unity.
Suppose that $B$ is defined in a neighborhood of $(0,0)$ as $f(x,y):=\sum_{i,j}c_{ij}x^iy^j=0$ in local coordinate. 
Since $B$ passes through $(0,0)$, $c_{00}=0$.
 By a direct computation,
 \begin{equation}
 \sigma^*f(x,y)=f_0(x,y)+f_1(x,y)\omega+f_2(x,y)\omega^2
 \end{equation}
 where 
 \begin{eqnarray*}
  f_0(x,y)&=&(c_{30}x^3+c_{12}xy)+\cdots\\
  f_1(x,y)&=&(c_{10}x)+(c_{40}x^4+c_{21}x^2y)+\cdots\\
 f_2(x,y)&=&(c_{20}x^2+c_{01}y)+(c_{50}x^5+c_{31}x^3y+c_{12}xy^2)+\cdots.
 \end{eqnarray*}
 As $B$ is invariant under the action of the cyclic group $\bZ_3$, it follows 
$B$ is defined by the vanishing of one of $f_0, f_1$ or $f_2$.  Since $B$ is smooth
at $(0,0)$, we conclude that $B$ is tangential to either $x+(c_{40}x^4+c_{21}x^2y)+\cdots=0$ 
or $y+c_{20}x^2+(c_{50}x^5+c_{31}x^3y+c_{12}xy^2)+\cdots=0$.  It means that if $B$ passes
through a point $P_i$, a fixed point of $\sigma$, 
$\hC$ is meeting precisely one of $E_{ij}, j=1,2$ on $\hX$.  After renaming, we may assume that
$\hC$ meets $E_{i1}$ but not $E_{i2}$.

It follows that in (9) $a_{i2}=0$ and $a_{i1}=-\hC\cdot E_{i1}/E_{i1}\cdot E_{i1}=\frac12$.

We compute 
\begin{eqnarray*}
K_{\hX}\cdot \hC&=&\tau^*K_X\cdot \hC=K_X\cdot C=1\\
\hC\cdot \hC&=&\tau^*C\cdot \tau^*C+\sum_{i=1}^2a_{i1}^2(E_{i1}\cdot E_{i1})=0.
\end{eqnarray*}
This implies that 
\begin{equation*}
2(g(\hC)-1)=\hC\cdot (K_{\hX}+\hC)=C\cdot K_X+\hC\cdot \hC=1>0,
\end{equation*}
contradicting $g(\hC)=g(C)=1$ as $C$ is a smooth elliptic curve.

\ms
\ni{\bf 5.4.} {\it Case (c). Fake projective plane as a non-regular covering of another complex $2$-ball quotient, possibly singular, of degree $3$:}

\ms
In this case, $M$ is a fake projective plane such that there is a non-regular degree $3$ covering
$p:M\rightarrow X$.  From the list given by the file \verb'registerofgps.txt' in the weblink of [CS], 
we check that there exist another fake projective plane $M'$ and a regular (normal) covering $p':M'\rightarrow X$ of degree $3$
for each pair of $M, X$ as above.
  Hence $X$ satisfies the properties listed in Lemma 6.

Assume again that $2K_M$ is not very ample so that a smooth curve $B$ of genus $3$ can be found  on $M$ as in {\it Case (b)}.  Since $M$ has Picard number $1$, the curve $p^{-1}(p(B))$
is connected.  Let $C=p(B)$.  
Let $\sigma:\tC\rightarrow C$ be the normalization of $C$.  The mapping $p$ can be factorized as $p=\sigma\circ q$, where $q:B\rightarrow \tC$.
Since $p:M\rightarrow X$ is a degree $3$ map, the degree of the mapping $q$ satisfies either, case (i), $\deg(q)=3$, or, case (ii), $\deg(q)=1$.  

Let us first rule out case (ii).

From Lemma 6, the singularities of $X$ consist of three
isolated singular points $Q_1, Q_2, Q_3$ of type $\frac13(1,2)$.   Consider now the non-regular mapping of degree
$3$ $p:M\rightarrow X$.  As $M$ is smooth, it follows from the singularity types of
$Q_i$, $i=1,2,3$ that $Q_i=p(P_i)$ for three points $P_1, P_2, P_3$ on $M$, and the local degree of $p$ around each $P_i$ is $3$.  

Since $p^{-1}(C)$ would have $2$ or $3$ irreducible components for which $B$ is one of them, we have to consider either case (iia), $p^{-1}(C)=B\cup B_2\cup B_3$ with irreducible $B_2, B_3$ or
case (iib), $p^{-1}(C)=B\cup B_4$ with irreducible $B_4$.  As each component is numerically equivalent to a multiple of $H$, they have non-trivial intersections, which leads to a singularity on $p^{-1}(C)$.

Assume first that $B$ intersects the set $\{P_i\}$.
Recall that $X=M'/\bZ_3$.  We also use the notation that $H_N$ denotes a generator of the Neron-Severi group modulo torsion in the case that $N$ has Picard number $1$, so that
a positive multiple of $H_N$ is an ample divisor.  Note that each of  $M, M'$ and $X$ has Picard number one.
In such a case, $H_{M'}$ descends to a $\bQ$-line bundle $H_X$ on $X$ with $H_X\cdot H_X=\frac13H_{M'}\cdot H_{M'}=\frac13.$
Hence $p^*H_X\equiv H_M$ on $M$.  Note that $K_X\equiv 3H_X$ and $K_M=3H_M$ as a line bundle.  It follows that 
$C\cdot H_X=p_*B\cdot H_X=B\cdot p^*H_X=1$ as $p|_B$ has degree $1$.  Hence $C\equiv K_X$.


Considering the Taylor expansion of a section of $H^0(X, 2K_X)$ at $Q_1$, three independent conditions are imposed for a non-trivial section to vanish at $Q_1$ to order at least $2$.  Since $h^0(X,2K_X)=4$ from Lemma 6,
we can find a section $s\in \Gamma(X,2K_X)$ such that
$s$ vanishes to order at least $2$ at $Q_1$.  Since $p^*K_X=K_M$, we know that $p^*s\in \Gamma(M,2K_M)$.  Since $p$ has local branching 
of order $3$ at $Q_1$, it follows that $p^*s$ vanishes to order at least $3\cdot 2=6$ at $P_1$.  It follows that $B$ intersects $p^*s$ to order at least $6$ at $P$.
Since $B\cdot 2K_M=6$, this implies that $B$ intersects $p^*s$ precisely to order $6$ at $P_1$ and they intersect nowhere else.  Moreover, it shows that $s$ vanishes at $Q_1$ precisely to order $2$, for
otherwise the intersection is larger than $6$.

On the other hand, $s\cdot C=2K_X\cdot K_X=6$.  Since $s$ and $C$ can intersect at most to order $2$ at $Q_1$, $s$ must intersect $C$ at some other point on $X$.  This implies that
$B$ intersects $p^*s$ at another point apart from $P_1$.  This contradicts the conclusion from the last
paragraph.  

Assume now that $B$ does not intersect the set $\{P_i\}$, so that $C$ does not intersect the set $\{Q_i\}$.  First we observe that $C$ cannot be smooth.  Since $p^{-1}(C)$ contains $B$ and $p:B\rightarrow C$ has degree $1$, there exists at most another irreducible component in $p^{-1}(C)$ apart from $B$.
However any other component is numerically an integral multiple of $H$ and hence has non-trivial intersection with $B$.  The intersection point gives rise to singularity of 
$p^{-1}(C)$.  Since $p$ is unramified outside the set $\{P_i\}$, $C$ must have singularity as well.  Since $q$ has degree one by our assumption, the singularity comes from the self-intersection of $\sigma(\tC)$.
Let $R$ be such a singular point, so that $p^{-1}(R)\cap B$ contains at least two points, say $S_1$ and $S_2$.  Again, we can find a section $s\in \Gamma(X,2K_X)$ vanishing at $R$ to order 
at least $2$.  We may assume that the vanishing order is $2$, since the proof for the case of order greater than $2$ is exactly the same.  In such a case,
$p^*s$ intersects $B$ at order at least $2$ at each of $S_1$ and $S_2$, since $p:M-\{P_i\}\rightarrow X-\{Q_i\}$ is unramified.  As $s\cdot C=6$ from earlier discussions, we know that $s$ intersect
$C$ at $4$ other points counted with multiplicity.  It follows that $p^*s$ intersects $B$ with intersection number at least $4$ at points other than $p^{-1}(R)\cap B$.  Together with the earlier
count, this implies that the $B\cdot p^*C\geqslant 8$, again contradict to the fact that $B\cdot p^*C=H\cdot 6H=6$.

Hence case (ii) is ruled out.


\ms
Consider now $\deg(q)=3$ as in case (i).  In such case, $C$ is linearly equivalent to
$H$ on $X$ since the Picard number of $X$ is $1$.  It follows from Hurwitz Formula that (8) still holds, except that the set of points $x$ such that $b_i(x)>0$ may or may not occur at the 
singularities of $X$.   Hence we conclude that $C$ is a smooth elliptic curve.  Recall that $p':M'\rightarrow X$ is a regular covering.  
In such a case, $B'=(p')^*C$ is linearly equivalent to $H$ on $M'$.  Using the result of Lemma 6 of [LY] as in the proof of (b) above, 
this implies that $B'$ is smooth of genus $3$.  Repeating the argument
of (b) to the regular covering $p':M'\rightarrow X$, we reach a contradiction as in (b).
Hence case (i) cannot happen as well.

\ms
In conclusion,  $B$ does not exist in our situation and hence $2K_M$ is very ample.

\ms
\ni{\bf 5.5.} {\it Case (d). Fake projective planes as a non-regular covering of another ball quotient of degree $21$:}

\ms
In this case, there is a non-regular degree $21$ covering
$p:M\rightarrow X$, and $M$ is not one of the cases in (b) or (c).   From the list given by the file \verb'registerofgps.txt' in the weblink of [CS], the only 
fake projective plane in this case is $M=(a=7, p=2,\emptyset,7_{21})$, which covers $X=(a=7, p=2,\emptyset)$.  
There exist another fake projective plane $M'=(a=7, p=2,\emptyset,D_3, 2_3)$ and a regular (normal) covering $p':M'\rightarrow X$ of degree $21$.

\begin{lemm}
The singularity set $\cS$ of $X$ consists of three points $Q_1, Q_2, Q_3$ of type $\frac13(1,2)$ and
one point $R$ of type $\frac17(1,3)$.  Let $\tau:\hX\rightarrow X$ be the minimal resolution of $X$.  Each $Q_i$ gives rise to a chain of two $-2$ curves $E_{ij}, j=1,2$ as in the case of (b).  The point $R$ 
is resolved to a chain of three rational curves $S_1, S_2, S_3$ of self-intersections $(-2)$, $(-2)$ and $(-3)$ respectively.  Moreover,
\begin{equation}
\tau^*K_X=K_{\hX}+\frac17S_1+\frac27S_2+\frac37S_3,
\end{equation}
$K_{\hX}^2=0$, and $c_2(\hX)=12$.
\end{lemm}

\ni{\bf Proof} From the existence of $p':M'\rightarrow X$, we know the structure of singular set as given in [K], see also [CS].  The others follow in a straight forward way, see also [K] or
[Y] for computations. 
Since $K_X$ is a $\bQ$ line bundle on $X$ with $K_X\cdot K_X=\frac3{7}$, we conclude that from the above that $K_{\hX}^2=0$, a fact found in [K].
Furthermore, equation (9) implies that $K_{\hX}\cdot\tau^*K_X=\frac37$.

\qed

Note that for each point $Q_i$, the set $p^{-1}(Q_i)$ consists of seven points $P_{ij}, j=1,\dots,7$, and $p^{-1}(R)$ consists of three points
$S_j, j=1,2,3$.

Let $C=p(B)$.  As the Picard number of $M$ is $1$, so is the one of $X$.   In this case, we know from [CS] that the lattice associated to $X$ can be lifted to $SU(2,1)$.  
Hence $H$ on $\tM$ descends to $H_X$ on $X$ and $H_M$ on $M$.  It follows that $p^*H_X=H_M$.  Since $p$ is unramified and
hence locally biholomorphic on $M-p^{-1}\cS$, we conclude that $C=p(B)$ is smooth on $X-\cS$. 
Hence $C\equiv kH$ for some positive integer $k$, where $H$ is a $\bQ$-line bundle with
$H\cdot H=\frac1{21}.$

Denote by $d$ the degree of $p|_B:B\rightarrow C$, which is the same as the degree of $B\rightarrow \tC$, the normalization of $C$ .  It follows that
$$1=B\cdot H_M=(p|_B)^*C\cdot H_M=dC\cdot H_X=dkH_X\cdot H_X=\frac {dk}{21}.$$
Hence 
$$dk=21.$$  
We conclude that $d$ can only take the value of $1$, $3$, $7$ or $21$.  In the following, we are going to eliminate these cases one by one.

In terms of notations in Lemma 7, we may write
\begin{equation}
\tau^*C=\hC+(\sum_{k=1}^3a_kS_k+\sum_{i=1}^3\sum_{j=1}^2b_{ij}E_{ij}),
\end{equation}
where $a_{k}$ and $b_{ij}$ are some integers.  By taking intersection with $K_{\hX}$ and making use of (9),
we get
\begin{eqnarray*}
K_{\hX}\cdot \tau^*C&=&K_{\hX}\cdot\hC+\sum_{k=1}^3a_kK_{\hX}\cdot S_k+\sum_{i=1}^3\sum_{j=1}^2b_{ij}K_{\hX}\cdot E_{ij}\\
&=&K_{\hX}\cdot\hC+\sum_{k=1}^3b_{ij}(\tau^*K_X-\frac17S_1-\frac27S_2-\frac37S_3)\cdot S_{k}\\
&=&K_{\hX}\cdot\hC+a_3,
\end{eqnarray*}
where we used the fact that $K_{\hX}\cdot E_{ij}=0$ since $E_{ij}$ is a $(-2)$ curve, $\tau^*K_X\cdot S_i=0=\tau^*K_X\cdot E_{ij}$, and the 
values of $S_i\cdot S_j$ given by Lemma 7.  Similarly
$$K_{\hX}\cdot \tau^*C= (\tau^*K_X-\frac17S_1-\frac27S_2-\frac37S_3)\cdot \tau^*C=3K_X\cdot C=\frac{k}{7}.$$
Hence
\begin{equation}
K_{\hX}\cdot\hC=\frac{k}7-a_3.
\end{equation}
From the above equation, we conclude that $\frac{k}7$ is an integer and hence $k$ can only be $7$ or $21$.

Assume first that $k=7$.  It follows that $d=3$.  From Riemann-Hurwitz Formula, we get as in (8) that
\begin{equation}
4=3\cdot2(g(\tC)-1)+\sum_{i=1}^lb_i,
\end{equation}
where $\tC$ is the normalization of $C$ and $l$ is the number of ramification points.   In this case, note that
$p^{-1}(C)$ has another component $B_1$ apart from $B$ since $p^{*}(C)\cdot p^{*}(C)=(7)^2H\cdot H>B\cdot B$.
In fact, $B_1\equiv 6H$ numerically.
Since $b_i$ can only take the values of $2$ or $6$ depending on the local branching order of $\tau|_{\hC}$, 
we conclude that $g(\tC)=1$, $l=2$ and $b_i=2$ for all $i$.  

Suppose that $a_3\neq0$.  It follows that $C$ passes through $R$.  As $R$ is a singularity of type $\frac17(1,3)$ and $M$
is smooth, a neighborhood of any point on $p^{-1}(R)$ gives rise to a local uniformization of $R$.  Since the degree of $p$ is $21$ and 
it is known that the ramification of $p$ consists only of points, cf. [CS], we conclude that $p^{-1}(R)$ consists of
three points $T_1, T_2, T_3$ and the local ramification order is $7$ at each $T_i$.   Hence $B$ has to pass one of the $T_i$'s
from our setting that $a_3>0$.  Without loss of generality, we may assume that $B$ passes through $T_1$.  Let $U_i$ be a small neighborhood of $T_i$.
Since $b_i=2$, it follows that $p|_{U_1}$ cannot be ramified over $R$,
which otherwise would give vanishing order of $7$.  Since we assume that $p^{-1}(C)=B_1\cup B$,
it follows that $B_1$ would intersect $B$ to order $6$ at one of the $S_i$.
As $B_1\cdot B=6H\cdot H=6$, it follows that $B_1$ does not intersect $B$ at any other point.  As $B$ is smooth on $M$
and $C=p(B)$ and $B_1$ does not intersect $B$ on $p^{-1}(X-\{R\})$, it follows that $C$ is smooth on $X-\{R\}$.
Since $B$ is smooth at $T_1$ and $p|_B:B\rightarrow C$ is unramified at $T_1$, we conclude that $C$ is smooth 
at $R$ as well.  Hence $C$ is smooth.  It follows that $p^{-1}(C)$ is also smooth.  In particular, $B_1$ is smooth except
at a set $F$ of self-intersection points of different local branches of $B_1$, where $F$ is necessarily a
subset of $p^{-1}(\cS)$.  Let $N_B$, $N_{B_1}$ be the normal bundle of $B$ and $B_1$ in $M$ respectively.  Denote by
$\Theta(N_B)$ the curvature with respect to the Poincar\'e metric.  Let $x\in B$ and $y\in p^{-1}(p(x))\in B_1$.
We equip the normal bundle $N_{B_1}(y)$ with a metric induced  from the metric of $N(B)(x)$, making use of the fact that $p$ is a local biholomorphism
at such points.  Hence the induced metric for the normal bundle of $B_1$ has curvature given by
$\Theta(B_1)(y)=\Theta(B)(x).$
It follows from the fact that $B$ is smooth that
$$\int_B\Theta(B)=H\cdot H=1.$$
Denote by $(B_1\cdot B_1)_o(y)$ the contribution to $B_1\cdot B_1$ given by the intersection of different local branches
of $B_1$ at $y\in F$.
Hence we get
\begin{equation}
B_1\cdot B_1=\int_{B_1-F}\Theta(B_1)+\sum_{y\in F}(B_1\cdot B_1)_o(y).
\end{equation}
It follows that $\int_{B_1-F}\Theta(B_1)=6\int_B\Theta(B)=6$,
and 
\begin{eqnarray*}
\sum_{y\in F}(B_1\cdot B_1)_o(y)&\geqslant& \sum_{i=1}^3(B_1\cdot B_1)_o(T_i)\\
&=&
\left(\begin{array}{c}6\\2\end{array}\right)+\left(\begin{array}{c}7\\2\end{array}\right)+\left(\begin{array}{c}7\\2\end{array}\right)\\
&=&57,
\end{eqnarray*}
where the second equality is computed from the self-intersection of $B_1$ at the three points $T_1, T_2$ and $T_3$ above $R$.
However the left hand side of (15) is $36$ since $B_1$ is numerically equivalent to $6H$.  We reach a contradiction.
Hence it cannot happen that $k=7$.

Suppose now that $a_3=0$.  As in the discussions in {\bf 5.4}, for each $i=1,\dots, l$, we may assume that $\hC$ intersects
$E_{i1}$ and does not intersect $E_{i2}$.  Hence $b_{i2}=0$ and $b_{i1}=\hC\cdot E_{i1}/E_{i1}\cdot E_{i1}=\frac12.$
Equation (13) now implies that $K_{\hX}\cdot\hC=1.$  It follows from (14) that 
either $l=5, g(\tC)=0$, or $l=2, g(\tC)=1$.  In either case, the right hand side of (14) is not an integer, which is a contradiction.


\ms
Assume now that $k=21$, we know that $d=1$.  The genus of the normalization $\tC$ of $C$ satisfies
$g(\tC)=g(B)=3$ and similarly the proper transform $\hC$ satisfies
$g(\widetilde{\hC})=g(\tC)=3$.
In this case, (13) is still valid with $k=21$, so that 
\begin{equation}
K_{\hX}\cdot\hC=\frac{21}7-a_3=3-a_3.
\end{equation}
From (12), we compute that
\begin{equation}
\hC\cdot \hC=\tau^*C\cdot\tau^*C+(\sum_{k=1}^3a_kS_k)^2+\sum_{i=1}^3(\sum_{j=1}^2b_{ij}E_{ij})^2.
\end{equation}

  From  Adjunction 
Formula,
\begin{eqnarray}
4&=&2(g(\widetilde{\hC})-1)=2(g(\hC)-\delta(\hC)-1)=\hC\cdot(K_{\hX}+\hC)-2\delta(\hC) \\
&=&4+(\sum_{k=1}^3a_kS_k)^2+\sum_{i=1}^3(\sum_{j=1}^2b_{ij}E_{ij})^2-a_3-2\delta(\hC)\nonumber
\end{eqnarray}
where we have used (16) and (17).  Here $a_3\geqslant 0$ and $\delta(\hC)\geqslant 0$,  
and $(\sum_{k=1}^3a_kS_k)^2+\sum_{i=1}^3(\sum_{j=1}^2b_{ij}E_{ij})^2$ is negative since intersection matrix associated 
to the contracted divisors is negative definite.  It follows from (18) that 
all these terms are $0$.  Hence $C$ is smooth and missed the singular set $\cS=\{R,Q_1, Q_2, Q_3\}$.  We know already that
$p:M\rightarrow X$ is of degree $21$ and unramified outside $p^{-1}(\cS)$.  Hence $p^{-1}(C)$ has at least another irreducible component $B_1$ apart from $B$.
$B_1$ must intersect $B$ since Picard number of $M$ is $1$.  This however is not possible since $p|_{M-\cS}$ is unramified and $C\cap \cS=\emptyset$.  The contradiction implies that $k$ cannot be $21$.

Since all the cases lead to contradiction, we conclude that such $B$ does not exist for $M=(a=7, p=2,\emptyset,7_{21})$.

\ms
\ni{\bf 5.6.} 
Theorem 3 is a consequence of the discussions for all the different cases from {\bf 5.2-5.5}.

\qed

\ms\ni
{\bf 5.7}
The above arguments cover all smooth surfaces of Euler number $3$ except for fake projective planes of minimal type, which consist of
 four lattices as mentioned in the Introduction.  We have the following result for such exceptional case.

\begin{lemm}
Let $\Gamma$ be a lattice associated to a fake projective plane $M$ of {\it minimal type}.  Then $2K_M$ is an embedding everywhere except that there may exist exactly one pair of points $p, q\in M$ so that $\Phi_{2K_M}(p)=\Phi_{2K_M}(q)$, where $p, q$ may be infinitesimally close.\\
\end{lemm}

\ms
\ni{\bf Proof}  Let $M$ be a general fake projective plane and $B$ be a section as given in Proposition 1(ii).
On $B$, consider long exact sequence 
\begin{eqnarray*}0&\rightarrow& H^0(B, 2K_M|_B-[p]-[q])\rightarrow H^0(B,2K_M|_B)\stackrel{\alpha}\rightarrow H^0(B, 2K_M|_p\oplus 2K_M|_q)\\&\rightarrow& H^1(B, 2K_M|_B-[p]-[q])\rightarrow
\end{eqnarray*}
From the Adjunction Formula, $K_B=K_M+B$.  Hence $2K_M|_B=K_B+K_M|_B-B|_B=K_B+2B+\tau|_\tau$.  Here and in the following,
we denote $B|_B$ and $\tau|_B$ simply by $B$ and $\tau$ when there is no danger of confusion.
Hence
$H^1(B, 2K_M|_B-[p]-[q])=H^1(B,K_B+2B-[p]-[q]+\tau)$, which is dual to $H^0(B, [p]+[q]-2B-\tau)$ from Serre Duality.  


There are two cases\\
Case (i),   $H^0(B, [p]+[q]-2B-\tau)=0$,\\
Case (ii),  $H^0(B, [p]+[q]-2B-\tau)\neq 0$.

Consider first Case (i).   In such case, $\alpha$ is surjective.
Consider now the long exact sequence
\begin{eqnarray*}
0&\rightarrow& H^0(M, 2K_M-B)\rightarrow H^0(M,2K_M)\stackrel{\beta}\rightarrow H^0(B, 2K_M|_B)\\
&\rightarrow& H^1(M, 2K_M-B)\rightarrow
\end{eqnarray*}
As $H^1(M, 2K_M-B)=0$ from Kodaira Vanishing Theorem, we conclude that $\beta$ is surjective.  The surjectivity of $\alpha$ and $\beta$
imply that $2K_M$ separates the points at $p$ and $q$, contradictory to our assumption. 
Hence $\Phi_{|2K_M|}$ gives an embedding

Consider now Case (ii).  Since $[p]+[q]-2B-\tau$ has degree $0$, Case (ii) can happen only if
$2B+\tau=[p]+[q]$ as an effective divisor.  Since $2B$ and $\tau$ are fixed, there can only be one pair of such $p$ and $q$.  This concludes
the proof of the lemma.

\qed
\ms




\ni{\bf 5.8.} \ni{\bf Proof of Theorem 2}  

The results of [Y4], see also  [Y5] or the revised version of [Y4] on the web for corrections, show that
a smooth compact complex $2$-ball quotient with $c_2=3$ has to be either a fake projective plane or a Cartwright-Steger surface.
The conclusion now follows from Theorem 3 and Lemma 6.

\qed


\ms
\ni{\bf 5.9.} \ni{\bf Proof of Theorem 1}  

From Lemma 1, we know that $c_1^2(M)$ is a positive multiple of $9$.    Suppose that $c_1^2(M)\geqslant 10$, it already
follows from the result of Reider [R] that $2K_M$ is very ample except for Case (i) described in Proposition 1.  That is, for two points $x_1, x_2\in M$ which
are not separated by $\Phi_{|2K_M|}$,
there is a curve $B$ passing through two points $x_1, x_2$ as described in Case (i).  In such case, again from Adjunction Formula,
$g(B)=2$.  The argument as in the second paragraph in the proof of Lemma 4 leads us
to a contradiction.  We conclude that $2K_M$ is very ample for $c_1^2(M)\geqslant 10$.  Theorem 1 now follows from Theorem 2 in the case of
$c_1^2(M)=9$.

\qed

\bs
\begin{center}
{\bf Appendix I} 
\end{center} 

In this appendix, we list the fake projective planes according to types discussed in the proof of Theorem 3.
The naming of the fake projective planes is given according to the naming of Cartwright-Steger in \verb'registerofgps.txt' in the weblink of [CS].
There are altogether $50$ lattices of $PU(2,1)$.  Further explanation of the nomenclature can be found in [PY] and [CS]. 

\begin{center}
{\bf Table (I): List of cases of fake projective planes}
\end{center}

{\scriptsize
$$\hskip 0in
\begin{array}{|c|c|}
\hline
M&Cases\\ \hline\hline
(a=1,p=5,\emptyset, D_3)& (b)\\ \hline
(a=1,p=5,\{2\},D_3)&(b)\\ \hline
(a=1,p=5,\{2I\})&(c)\\ \hline
(a=2,p=3,\emptyset,D_3)&(b)\\ \hline
(a=2,p=3,\{2\},D_3))&(b)\\ \hline
(a=2,p=3,\{2I\})&(c)\\ \hline
(a=7,p=2,\emptyset, D_32_7)&(b)\\ \hline
(a=7,p=2,\emptyset,7_{21})&(d)\\ \hline
(a=7,p=2,\emptyset,D_3X_7)&(b)\\ \hline
(a=7,p=2,\{7\},D_32_7)&(b)\\ \hline
(a=7,p=2,\{7\},D_37_7)&(b)\\ \hline
(a=7,p=2,\{7\},D_37'_7)&(b)\\ \hline
(a=7,p=2,\{7\},7_{21})&(c)\\ \hline
(a=7,p=2,\{3\},D_3)&(b)\\ \hline
(a=7,p=2,\{3\},3_3)&(c)\\ \hline
(a=7,p=2,\{3,7\},D_3)&(b)\\ \hline
(a=7,p=2,\{3,7\},3_3)&(c)\\ \hline
(a=7,p=2,\{5\})&\mbox{\mbox{min type}}\\ \hline
(a=7,p=2,\{5,7\})&\mbox{\mbox{min type}}\\ \hline
(a=15,p=2,\emptyset,D_3)&(b)\\  \hline
(a=15,p=2,\emptyset,3_3)&(c)\\ \hline
(a=15,p=2,\{3\},D_3)&(b)\\ \hline
(a=15,p=2,\{3\},3_3)&(b)\\ \hline
(a=15,p=2,\{3\},(D3)_3)&(b)\\ \hline
(a=15,p=2,\{5\},D_3)&(b)\\ \hline
(a=15,p=2,\{5\},3_3)&(c)\\ \hline
(a=15,p=2,\{3,5\},D_3)&(b)\\ \hline
(a=15,p=2,\{3,5\},3_3)&(b)\\ \hline
(a=15,p=2,\{3,5\},(D3)_3)&(b)\\ \hline
(a=23,p=2,\emptyset)&\mbox{\mbox{min type}}\\ \hline
(a=23,p=2,\{23\})&\mbox{\mbox{min type}}\\ \hline
(\cC_2,p=2,\emptyset, d_3D_3)&(b)\\ \hline
(\cC_2,p=2,\emptyset, D_3X_3)&(b)\\ \hline
(\cC_2,p=2,\emptyset, (dD)_3 X_3)&(b)\\ \hline
(\cC_2,p=2,\emptyset, (d^2D)_3X_3)&(b)\\ \hline
(\cC_2,p=2,\emptyset, d_3X'_3)&(c)\\ \hline
(\cC_2,p=2,\emptyset, X_9)&(c)\\ \hline
(\cC_2,p=2,\{3\},d_3D_3)&(b)\\ \hline
(\cC_{10},p=2,\emptyset,D_3)&(b)\\ \hline
(\cC_{10},p=2,\{17-\},D_3)&(b)\\ \hline
(\cC_{18},p=3,\emptyset,d_3D_3)&(b)\\ \hline
(\cC_{18},p=3,\{2\},D_3)&(b)\\ \hline
(\cC_{18},p=3,\{2\},(dD)_3)&(b)\\ \hline
(\cC_{18},p=3,\{2\},(d^2D)_3)&(b)\\ \hline
(\cC_{18},p=3,\{2I\})&(c)\\ \hline
(\cC_{20},\{v_2\},\emptyset,D_32_7)&(b)\\ \hline
(\cC_{20},\{v_2\},\{3+\},D_3)&(b)\\ \hline
(\cC_{20},\{v_2\},\{3+\},\{3+\}_3)&(c)\\ \hline
\cC_{20},\{v_2\},\{3-\},D_3)&(b)\\ \hline
(\cC_{20},\{v_2\},\{3-\},\{3-\}_3)&(c)\\ \hline
\end{array}
$$
}

\bs
\begin{center}
{\bf Appendix II} 
\end{center} 

In the proof of Theorem 3 for Cartwright-Steger surface in {\bf 5.2}, the details for elimination of curve $B$ satisfying (i) of Proposition 1 was not given, as cautioned by Jonghae Keum in an email.  The details for non-existence of $B$ is given below,
consisting of two email messages from the author to Keum in December 2017 except the last paragraph.  The author is grateful to Keum for the correspondence.

Assume for the sake of proof by contradiction that $B$ as above exists on Cartwright-Steger surface $M$, so that
 \begin{equation}
0=B\cdot B=0,\
 2=K\cdot B
 \end{equation}
 In such case, $B$ has genus $2$, and is smooth and irreducible from the fact
 that  $M$ is complex hyperbolic.  $h^0(M,B)=1$, for otherwise the ratio of two linearly independent ones
 would give a fibration from $B\cdot B=0$, and the argument in the paper leads to a singular fiber with normalization
 of genus $\leqslant 1$, contradicting hyperbolicity.  
 
 Let $F$ be the divisor class representing fiber of the Albanese map.  From [CKY],
\begin{equation}
F\equiv -E_1+5E_2.
\end{equation}
Since $B$ is effective, we also have
\begin{equation} 
F\cdot B=n
\end{equation}
for some non-negative integer $n$.  Here we use Zariski Lemma if $B$ is a component of $F$.

From [CKY], there are five explicit divisors $E_1, E_2, E_3$ with $E_3\equiv K$ and $E_1+E_2\equiv 2K$, and $C_1, C_2$.
 We have the following intersection matrix for the curves $E_1$, $E_2$, $E_3$ and $C_1, C_2$:
\begin{equation}\label{eq:intersectionmatrix}
\begin{pmatrix}
5&13&9&11&11\\
13&5&9&7&7\\
9&9&9&9&9\\
11&7&9&-1&17\\11&7&9&17&-1
\end{pmatrix}.
\end{equation}
From the above, we can see easily that the following form an orthogonal basis of $NS/\tor$ with respect to intersection pairing,
$$K, E_1-E_2, D:=C-K+\frac14(E_1-E_2).$$
We have
\begin{equation}
K\cdot K=9, \ (E_1-E_2)\cdot (E_1-E_2)=-16, D\cdot D=-9.
\end{equation}
From (20), we can write
\begin{equation}
F\equiv 4K-3(E_1-E_2).
\end{equation}
In terms of the basis chosen, we write 
\begin{equation}
B\equiv aK+b(E_1-E_2)+cD
\end{equation}
for some rational numbers $a, b, c$.

From (19) and (25), $a=\frac29$.   

From (19) and (23),  we get $0=\frac49-16 b^2-9c^2$, which can be rewritten as 
\begin{equation}
(6b)^2+(\frac92 c)^2=1.
\end{equation}
From (21), (24) and (25), 
\begin{equation}
6b=\frac n8-1
\end{equation}
Hence we conclude that
\begin{equation}
c=\pm\frac1{36}\sqrt{(16-n)n}
\end{equation}
Since $c$ is rational, the only possibilities are\\
(i) $n=0$, \\
(ii) $n=16$ or\\
(iii) $n=8$.

Consider first Case (i) that $n=0$.  In such case, $F$ is represented by a general fiber of the Albanese fibration $\alpha$, we see that $B$ lies in the fiber.
 Hence $B$ has to  be an irreducible component of a fiber $F_o$ of $\alpha$. $B$ cannot be a fiber from genus consideration.
From [CKY],  any fiber of $\alpha$ is either stable or irreducible.  Hence the only possibility is that $F_o$ is connected and has more than one irreducible 
component.  In such case, it leads to contradiction since (19) leads to $B\cdot F_o=B\cdot F_1>0$, where $F_1$ is the union of the other components of $F_o$
apart from $B$.

Alternatively, in this case we know from (27), (28) that $b=-\frac16$ and $c=0$.  Hence $B\equiv\frac29 K-\frac16(E_1-E_2)$.  It follows that
$$B\cdot C=\frac29(9)-\frac16(11-7)=2-\frac23,$$
which is not an integer, a contradiction.

For Case (ii), we know from (27), (28) that $b=\frac16$ and $c=0$.  Hence $B\equiv\frac29 K+\frac16(E_1-E_2)$.  It follows that
$$B\cdot C=\frac29(9)+\frac16(11-7)=2+\frac23,$$
which is not an integer, a contradiction.
 
Consider now Case (iii).  In such case, $n=8$, $b=0$ and $c=\pm\frac29$ from (27) and (28).  Hence
\begin{equation}
B\equiv\frac29(K\pm(C-K+\frac14(E_1-E_2))).
\end{equation}

Now if we choose at the beginning $C=C_1$, it follows  that 
\begin{equation}
B\equiv\frac29(C_1+\frac14(E_1-E_2)), \  (\mbox{resp.} \ B\equiv\frac29(2K-C_1-\frac14(E_1-E_2)) )
\end{equation}
It follows from (23) that
$B\cdot C_1=0$ (resp. $B\cdot C_2=0$).  
In conclusion, $B$ as given in Proposition 1(i) has to satisfy for either $C=C_1$ or $C_2$ that
\begin{equation}
B\cdot C=0
\end{equation}

As $h^0(M,B)=1$, we have either \\
Case (a) $B$ is invariant under $C_3$, the automorphism
 of $M$, or \\
Case (b) $\sigma B=B+\tau$ for some $3$-torsion line bundle $\tau$, where $\langle\sigma\rangle=C_3$.
 
Consider first Case (a).  As the three explicit divisors $E_1, E_2, E_3$ are each invariant under $C_3$, the set $B\cap E_3$ is invariant under $C_3$.  As $B\cdot E_3=2$,
 the intersection occurs either at one point with multiplicity $2$, or two points with multiplicity $1$, and the points
 are fixed by $C_3$.
  In the notation of [CKY], Prop 5.5, the fixed points of $C_3$ contained in
 $E_i$'s are denoted by
  $O_1, O_2, O_3$.
 From [CKY], the number of local branches of
 $E_1, E_2$ and $E_3$ at the points $O_1, O_2, O_3$ are $(3,1,2)$, $(2,1,3)$ and
 $(1,4,1)$ respectively.   Hence the only possibility is $B\cap E_3=\{O_1, O_2\}$.  Similar argument
 for $B\cap(E_1+E_2)$ implies that the intersection contains $1$ or $4$ fixed points under the action of $C_3$.
 Again, as $O_1, O_2, O_3$ are the only candidates, it can only happen that the intersection occurs only at
 one point, with multiplicity between $1$ and $4$.  From the number of local branches mentioned above, it can
 only occur at $O_3$.  Altogether this means that that $O_1, O_2, O_3\in B$.  However, the intersection of $B$ with
 $2E_3$ is then at least $2(1+4+1)$, contradicting $B\cdot 2E_3=4$ and eliminating Case (a).

Consider now case (b).  There are  three disjoint effective divisors $B_1, B_2, B_3$ corresponding to sections
 $s_1\in\Gamma(M,B), s_2\in\Gamma(M,B+\tau), s_3\in\Gamma(M,B+2\tau)$ for a $3$-torsion $\tau$.  It follows that $s_1^2$ and $s_2s_3$ are
 two linearly independent sections of $\Gamma(M,2B)$.  We now have a fibration $\pi:M\rightarrow \bP^1$ given
 by $as_1^2+bs_2s_3$, where all fibers are disjoint with trivial self intersection, and are sections of $\Gamma(M,2B)$ with arithmetic genus $3$
 given by Adjunction Formula.  
 Let $x\in C$.  There exists a fiber $F_s$ of $\pi$ containing $x$.  As $F_s\cdot C=2B\cdot C=0$, 
 we conclude that $C$ is a component of $F_s$.  This however contradicts the fact that $F_s$ has arithmetic genus $3$, while
 $C$ has arithmetic genus $5$ by Adjunction Formula and (22). 
 Hence Case (b) is eliminated.  In conclusion, $B$ as in Proposition (i) does not exist.
 
\bigskip
\noindent{\bf References} 

\bs
\ni [BHPV] Barth, W. P., Hulek, K., Peters, C. A. M., Van de
Ven, A., Compact complex surfaces. Second edition. Ergebnisse der
Mathematik und ihrer Grenzgebiete. 3. Folge. A Series of Modern
Surveys in Mathematics 4. Springer-Verlag, Berlin, 2004.

\ms
\ni [BCP] Bauer, I., Catanese, F., Pignatelli, R., Surfaces of general type with geometric genus zero: a survey. Complex and differential geometry, 1-48, Springer Proc. Math., 8, Springer, Heidelberg, 2011. 

\ms
\ni [Bom]
Bombieri, E., Canonical models of surfaces of general type. Publ. Math., Inst. Hautes \'Etud. Sci. 42, 171-220 (1972)

\ms
\ni [Bor] Borrelli, G., The classification of surfaces of general type with nonbirational bicanonical map. J. Alg. Geo. 16 (2007), no. 4, 625-669.

\ms
\ni [CKY] Cartwright, D., Koziarz, V., Yeung, S.-K, On the Cartwright-Steger surface, to appear in Jour. Alg. Geo., arXiv:1412.4137.

\ms
\ni [CS] Cartwright, D., Steger, T., Enumeration of the $50$ fake projective planes, C. R. Acad. Sci. Paris, Ser. 1,
348 (2010), 11-13, see also \\
http://www.maths.usyd.edu.au/u/donaldc/fakeprojectiveplanes/

\ms
\ni
[DM] Deligne, P., Mostow, G. D.,  Monodromy of hypergeometric functions and non-lattice
integral monodromy,  Publications Mathematiques. Institut de Hautes Etudes Scientifiques
63 (1986): 5-89.

\ms
\ni [DD] Di Brino, G., Di Cerbo, L., Exceptional collections and the bicanonical map of Keums fake projective planes, arXiv:1603.04378v1.

\ms
\ni [HT] Hwang, J.-M., To, W.-K., On Seshadri constants of canonical bundles of compact complex hyperbolic spaces, Comp. Math. 118 (1999), 203-215.

\ms
\ni [K] Keum J., Quotients of fake projective planes. Geom Topol, 12 (2008), 2497-2515.

\ms
\ni
[LY] Lai, C. L., Yeung, S.-K., Exceptional collections on some fake projective planes, submitted,
http://www.math.purdue.edu/\%7Eyeung/papers/FPP-van\_r.pdf

\ms
\ni [Le] Le Vavasseur, R., Sur le syst\'eme d'\'equations aux d\'eriv\'e es partielles simultan\'e es auxquelles
satisfait la s\'erie hyperg\'eom\'etieque \`a deux variables, J. Fac. Sci. Toulouse VII (1896), 1-205.

\ms
\ni [Li]  Liu, K., Geometric height inequalities. Math. Res. Lett. 3 (1996), 693-702.

\ms
\ni [MP] Mendes Lopes, M., Pardini, R., The bicanonical map of surfaces with $p_g = 0$ and
$K^2 \geqslant 7$. II. Bull. London Math. Soc. 35, (2003), 337-343.

\ms
\ni
[M] Mumford, D.,  An algebraic surface with $K$ ample, $K^2=9$, $p_g=q=0.$
 Amer. J.
Math. 101(1979), 233--244.

\ms
\ni [P] Picard, E., Sur les fonctions hyperfuchsiennes provenant des s\'eries hyperg\'eom\'etriques de
deux variables,  Ann. ENS III, 2 (1885), 357-384.

\ms
\ni [PY] Prasad, G., Yeung, S.-K.,  Fake projective planes, Inv. Math. 168 (2007), 321--370; Addendum, ibid 182 (2010), 213-227.

\ms
\ni [R] 
Reider, I.,  Vector bundles of rank 2 and linear systems on algebraic surfaces. Ann. of Math. (2) 127 (1988), 309-316.

\ms
\ni [Y1] Yeung, S.-K., Very ampleness of line bundles and canonical embedding of coverings of manifolds. Comp. Math. 123 (2000), 209-223. 

\ms
\ni [Y2]  Yeung, S.-K., Classification and construction of fake projective planes, Handbook of geometric analysis, No. 2, 391-431, Adv. Lect. Math. (ALM), 13, Int. Press, Somerville, MA, 2010.

\ms
\ni [Y3]  Yeung, S.-K., Exotic structures arising from fake projective planes. Sci. China Math. 56 (2013), no. 1, 43-54, addendum, ibid, 58 (2015), no. 11, 2473-2476.

\ms
\ni [Y4] Yeung, S.-K., Classification of surfaces of general type with Euler number $3$, J. Reine Angew. Math. 697 (2014), 1-14. 
corrected version,\\ http://www.math.purdue.edu/\%7Eyeung/papers/fake-euler.pdf

\ms
\ni [Y5] Yeung, S.-K., Foliations associated to harmonic maps and some complex two-ball quotients, Sci. China Math., 60(2017), 1137-1148.

\end{document}